\documentclass[12pt]{amsproc}
\usepackage{amsmath,amssymb,amscd}

\begin{document}

\title[A Remark About Supramenability and the Macaev Norm]{A Remark About Supramenability and the Macaev Norm}
\author{Dan-Virgil Voiculescu}
\address{D.V. Voiculescu \\ Department of Mathematics \\ University of California at Berkeley \\ Berkeley, CA\ \ 94720-3840}
\thanks{Research supported in part by NSF Grant DMS-1301727.}
\date{}

\begin{abstract}
We show that a finitely generated group $G$ which satisfies a certain condition with respect to the Macaev norm is supramenable. The condition is equivalent to the existence of quasicentral approximate unit with respect to the Macaev norm relative to the left-regular representation of the group and has been studied by the author in connection with perturbation questions for Hilbert space operators. The condition can be also viewed as an analogue with respect to the Macaev norm of Yamasaki parabolicity. We also show that existence of quasicentral approximate units relative to the Macaev norm for $n$-tuples of operators is not preserved when taking tensor products.
\end{abstract}

\maketitle

\section{Introduction}
\label{sec1}

The Macaev norm is also known as the Lorentz $(\infty,1)$-norm, the name of Macaev being associated with the corresponding ideal of operators \cite{5}. In this note we exhibit a condition involving the Macaev norm which implies the supramenability (\cite{3}, \cite{11}) of a finitely generated group $G$. The Macaev ideal plays an important role in certain operator theory questions \cite{5} and in non-commutative geometry \cite{4}, where its dual is the habitat of the Dixmier trace.

The condition for supramenability amounts to the vanishing of an invariant $k_{\infty}^-$ which appeared in our work (\cite{13}, \cite{14}, \cite{15}, \cite{16}, \cite{17}) on Hilbert space operators modulo normed ideals of compact operators. In the case of a finitely generated group $G$ the condition can also be stated in terms of functions on $G$, that is without dealing with operators. It can also be viewed as an analogue with the Macaev norm replacing the $p$-norm, of Yamasaki $p$-parabolicity \cite{18}. Note also that in \cite{15}, \cite{17} we had found that the non-vanishing of $k_{\infty}^-$ is connected to the entropies of Kolmogorov--Sinai and Avez \cite{1}, and R.~Okayasu in \cite{8}, \cite{9}, \cite{10} further explored connections with entropy, for connections of Avez-entropy to operator algebras, see also \cite{2}. In particular on the groups for which our condition implies supramenability the finitary random walks have zero entropy, which has implications for the Poisson boundary \cite{6}.

Supramenability is a strong requirement about the existence of finitely additive invariant measures introduced by J.~M. Rosenblatt \cite{11}. We make essential use of a recent result of Kellerhals--Monod--R$\phi$rdam \cite{7} which makes precise a free semigroup obstruction to supramenability.

After the introduction in Section~\ref{sec2} we prove a version of the main result. In Section~\ref{sec3} we then explain various equivalent forms of the result. In a remark at the end, we bring up the question of the condition being possibly equivalent to supramenability. In Section~\ref{sec4} we give an example of $n$-tuples of operators for which the existence of quasicentral approximate units relative to the Macaev ideal does not carry over to their tensor product.

\section{A version of the main result}
\label{sec2}

In this section we prove a version of the main result and leave the reformulations, which will provide a better perspective, for the next section. We also give a generalization of a part of the proof.

By $(\ell_{\infty}^-(X),|\ |_{\infty}^-)$ we denote the Banach space of functions $f: X \to {\mathbb R}$ vanishing at infinity endowed with the norm
\[
|f|_{\infty}^- = \sum_{k \in {\mathbb N}} f^*(k)k^{-1}
\]
where $f^*(1) \ge f^*(2) \ge \dots$ is the decreasing rearrangement of $|f|$. Here $|\ |_{\infty}^-$ is the Macaev or Lorentz $(\infty,1)$ norm. By ${\mathcal F}(X)$ we denote the real-valued functions on $X$ with finite support.

Throughout $G$ will be a finitely generated group and $K \subset G$ some finite generator.

\bigskip
\noindent
{\bf Theorem 2.1.} {\em 
If $G$ is not supramenable then
\[
\inf\left\{ \max_{g \in K} |\alpha(g)f-f|_{\infty}^- \mid f \in {\mathcal F}(G),\ f(e) = 1\right\} > 0
\]
where $(\alpha(g)f)(g_1) = f(g_1g)$.
}

\bigskip
\noindent
{\bf {\em Proof.}} Assume $G$ is not supramenable and that
\[
\inf\left\{ \max_{g \in K} |\alpha(g)f-f|_{\infty}^- \mid f \in {\mathcal F}(G),\ f(e) = 1\right\} = 0.
\]
By Theorem~1.1 of \cite{7} the failure of supramenability implies the existence of a Lipschitz injective map $\rho: {\mathbb F}_2 \to G$, where ${\mathbb F}_2$ is the free group on $2$ generators. Adjusting $\rho$ with a right shift on $G$, we may assume $\rho(e) = e$. On the other hand there is a sequence $f_n \in {\mathcal F}(G)$, $f_n(e) = 1$ so that
\[
\lim_{n \to \infty} \max_{g \in K} |\alpha(g)f_n - f_n|_{\infty}^- = 0.
\]
Since $\rho$ is Lipschitz, if $g_1,g_2$ are the generators of ${\mathbb F}_2$ and $K = K^{-1}$, there is a number $M \in {\mathbb N}$ so that
\[
\rho(g_ph) \in K^M\rho(h)
\]
for $p \in \{1,2\}$, $h \in {\mathbb F}_2$. If $f \in {\mathcal F}(G)$ we have
\[
\begin{aligned}
|f(\rho(h)) - f(\rho(g_ph))| &\le \max_{t \in K^M} |f(\rho(h)) - f(t\rho(h))| \\
&\le \left(\left(\sum_{t \in K^M} |\alpha(t)f-f|\right) \circ \rho \right)(h).
\end{aligned}
\]
Since $|g \circ \rho|_{\infty}^- \le |g|_{\infty}^-$ and
\[
\begin{aligned}
\left| |\alpha(t)f-f|\right|_{\infty}^- &= |\alpha(t)f-f|_{\infty}^- \\
&\le M \max_{k \in K} |\alpha(k)f-f|_{\infty}^-
\end{aligned}
\]
when $t \in K^M$, we infer that
\[
\begin{aligned}
|\alpha(g_p)(f \circ \rho) - f\circ \rho|_{\infty}^- &\le \left|\left( \sum_{t \in K^M} |\alpha(t)f - f|\right) \circ \rho \right|_{\infty}^- \\
&\le \sum_{t \in K^M} |\alpha(t)f - f|_{\infty}^- \\
&\le M|K|^M \max_{k \in K} |\alpha(k)f-f|_{\infty}^-.
\end{aligned}
\]
This implies
\[
(f_n \circ \rho)(e) = 1,\ n \in {\mathbb N}
\]
and
\[
\lim_{n \to \infty} \max_{p \in \{1,2\}} |\alpha(g_p)(f_n \circ \rho)-f_n \circ \rho|_{\infty}^- = 0.
\]
Thus we would have
\[
\inf\left\{ \max_{p=1,2} |\alpha(g_p)f-f|_{\infty}^- \mid f \in {\mathcal F}({\mathbb F}_2),f(e)=1\right\} = 0.
\]

That this is false, can be seen using equivalent rephrasings either from \cite{14} or Corollary~6 of \cite{17}. For the reader's convenience we sketch here the argument for this, so that the proof is self-contained. On ${\mathbb F}_2$ consider the functions $H_p$, $p = 1,2$, which are $\ne 0$ on $g_pS$, where $S$ is the semigroup generated by $g_1,g_2,e$ and $H_p(g_{i_1}\dots g_{i_k}) = \delta_{p,i_1}2^{-k}$. Then
\[
(\alpha(g_1^{-1})H_1 - H_1) + (\alpha(g_2^{-1})H_2 - H_2)
\]
is the indicator function $\chi_{\{e\}}$ of the set $\{e\}$. The functions $H_1,H_2$ are in the dual of $\ell_{\infty}^-$ which is $\ell_1^+$ with norm 
\[
\sup_{k \in {\mathbb N}} \left(\left( \sum_{1 \le j \le k} f^*(j)\right)\left( \sum_{i \le j \le k} j^{-1}\right)^{-1}\right)\,.
\]
Then
\[
\begin{aligned}
\sum_{p = 1,2} |\alpha(g_p)f-f|_{\infty}^- &\ge C \left| \sum_{\underset{p = 1,2}{g \in {\mathbb F}_2}} (\alpha(g_p)f-f)(g)H_p(g)\right| \\
&= C \left| \sum_{\underset{p = 1,2}{g \in {\mathbb F}_2}} (f(g)(\alpha(g_p^{-1})H_p-H_p)(g)\right| \\
&= C|f(e)| = C
\end{aligned}
\]
where $C > 0$. Thus if $f_n \in {\mathcal F}({\mathbb F}_2)$, $f_n(e) = 1$ we cannot have
\[
\lim_{n \to \infty} \max_{p = 1,2}|\alpha(g_p)f_n-f_n|^-_{\infty} = 0.
\]
\qed

\bigskip
\noindent
{\bf Corollary 2.1.} {\em 
If
\[
\inf\left\{ \max_{k \in K}|\alpha(k)f-f|_{\infty}^-| \mid f \in {\mathcal F}(G),f(e) = 1\right\} = 0
\]
then $G$ is supramenable.
}

\bigskip
The part of the proof involving the use of the Lipschitz property is seen immediately to work also for the rearrangement invariant norms and will be recorded as the second theorem of this section.

Let $\Phi(\lambda_1 \ge \lambda_2 \ge \dots)$ where $\lambda_j \ge 0$, with only finitely many $\ne 0$, be a norming function (see \cite{5} or \cite{12}) and let $|\ |_{\Phi}$ the norm on ${\mathcal F}(X)$
\[
|f|_{\Phi} = \Phi(f^*(1) \ge f^*(2) \ge \dots)
\]
where $f^*(1) \ge f^*(2) \ge \dots$ is the decreasing rearrangement of $(|f(x)|)_{x \in X}$.

\bigskip
\noindent
{\bf Theorem 2.2.} {\em 
Let $G_1,G_2$ be groups with finite generating sets $K_1,K_2$ and assume there is an injective Lipschitz map $\rho: G_1 \to G_2$. Then
\[
\inf\left\{ \max_{g \in K_1} |\alpha(g)f-f|_{\Phi} \mid f \in {\mathcal F}(G_1),f(e) = 1\right\} > 0
\]
implies
\[
\inf\left\{ \max_{g \in K_2} |\alpha(g)f-f|_{\Phi} \mid f \in {\mathcal F}(G_2),f(e) = 1\right\} > 0.
\]
}

\section{Perturbation theory and reformulations}
\label{sec3}

The condition involving the Macaev norm also appeared in our note \cite{17} in connection with entropy of random walks and we explained there how it relates to other equivalent conditions in view of (\cite{14}, \cite{16}). For the reader's convenience we briefly state these equivalences again here, but leaving out details which can be found in (\cite{14}, \cite{16} or \cite{17}).

If $\Phi(\lambda_1\ge \lambda_2 \ge \dots)$ is a norming function and $X$ a finite rank operator on a complex separable infinite dimensional Hilbert space ${\mathcal H}$ the norm $|\ |_{\Phi}$ is defined by $|X|_{\Phi} = \Phi(s_1(X),s_2(X),\dots)$ where $s_j(X)$ are the eigenvalues of $(X^*X)^{1/2}$ in decreasing order and with repetitions according to multiplicity. If $\tau = (T_1,\dots,T_n)$ is an $n$-tuple of bounded operators on ${\mathcal H}$ the invariant $k_{\Phi}(\tau)$ is defined by
\[
k_{\Phi}(\tau) = \liminf_{A \in R_1^+({\mathcal H})} \max_{1 \le j \le n} |[T_j,A]|_{\Phi}
\]
where $R_1^+({\mathcal H})$ denotes the finite-rank positive contraction operators on ${\mathcal H}$ endowed with its natural order (\cite{13}, \cite{14}). The number $k_{\Phi}(\tau)$ is a kind of ``size ${\mathcal G}_{\Phi}$ measure'' of $\tau$, where ${\mathcal G}_{\Phi}$ is the normed ideal associated with $\Phi$ (\cite{5}, \cite{12}). For the norming function corresponding to the Macaev norm $|\ |_{\infty}^-$ the number will be denoted by $k_{\infty}^-(\tau)$. Among the remarkable properties of $k_{\infty}^-(\tau)$ (\cite{14}, \cite{15}, \cite{16}, \cite{17}) is the fact that if ${\mathcal G}_{\Phi} \supset {\mathcal C}_{\infty}^-$ and ${\mathcal G}_{\Phi} \ne C_{\infty}^-$ then $k_{\Phi}(\tau) = 0$ for all $\tau$. Note that $k_{\infty}^-(\tau)$ is always finite and that it is non-zero, for instance if $\tau$ is an $n$-tuple of Cuntz isometries, $n \ge 2$.

In \cite{15}, \cite{17} we found connections between $k_{\infty}^-$ and the Kolmogorov--Sinai entropy and the Avez-entropy. This direction was further explored in the work of R.~Okayasu (\cite{8}, \cite{9}, \cite{10}).

Let $G$ be a group with a finite generator $K \subset G$ and let $\lambda$ be the left regular representation of $G$ on $\ell^2(G)$. Then the fact that for a given $\Phi$ we have $k_{\Phi}(\lambda(K)) > 0$ or $k_{\Phi}(\lambda(K)) = 0$ is independent of the choice of $K$, that is it is a property of $G$.

The following conditions on $G$ are equivalent:

\begin{itemize}
\item[k.1.] $k_{\Phi}(\lambda(K)) > 0$
\item[k.2.] $\displaystyle{\inf\left\{ \max_{g \in K} |f-\alpha(g)f|_{\Phi} \mid f \in {\mathcal F}(G),f(e)=1\right\} > 0}$
\item[k.3.] $\displaystyle{\inf \left\{ \max_{g \in K} |f-\alpha(g)f|_{\Phi} \mid f \in {\mathcal F}(G),0 \le f \le 1,f(e)=1\right\} > 0}$
\item[k.4.] there is a constant $C > 0$ so that
\[
\|f\|_{\infty} \le C \max_{y \in K} |f-\alpha(g)f|_{\Phi}
\]
if $f \in {\mathcal F}(G)$.
\item[k.5.] if $\Phi^d$ is the dual norming function of $\Phi$ (see \cite{4}, \cite{8}) then there are $f_k \in \ell_{\Phi^d}(G)$ (which is the dual of $\ell_{\Phi}(G)$) so that $\displaystyle{\sum_{k \in K} (\alpha(k^{-1})f_k-f_k)(g) = \delta_{g,e}}$.
\end{itemize}

Condition k.4 in the case where $\Phi$ is the $p$-norm is Yamasaki's hyperbolicity condition \cite{18}.

For instance Theorem~2.2 can be reformulated as saying that if there is a Lipschitz injection $\rho: G_1 \to G_2$ then
\[
k_{\Phi}(\lambda(K_1)) > 0 \Rightarrow k_{\Phi}(\lambda(K_2)) > 0,
\]
where $K_1,K_2$ are the generating sets of $G_1$ and $G_2$, respectively.

\bigskip
\noindent
{\bf Remark 3.1.} It is natural to ask whether the supramenability of $G$ is actually equivalent to $k_{\infty}^-(\lambda(K)) = 0$. Note that it is possible that one, or both, of these two properties be also equivalent to some subexponential growth condition on $G$. Since in \cite{17} we showed that $k_{\infty}^-(\lambda(K)) > 0$ if there is a finitely supported probability measure $\mu$ on $G$ for which the Avez-entropy of the random walk $h(G,\mu) > 0$, the equivalence would fail if there is a supramenable $G$ for which $h(G,\mu) > 0$ for some $\mu$. On the other hand, the counterexample in the next section was inspired by \cite{3}, a transfer from supramenable to $k^-_{\infty}$.

\section{Tensor product of $n$-tuples of operators}
\label{sec4}

In this section we show that $k_{\infty}^-(\tau) = k_{\infty}^-(\tau') = 0$ does not imply $k_{\infty}^-(\tau \otimes \tau') = 0$. Actually the assumption can be strengthened, as will be seen below. The counterexample is inspired by an example of L.~Bartholdi (74 Example in \cite{3}) of two supramenable metric spaces, the product of which fails to be supramenable.

\bigskip
\noindent
{\bf Proposition 4.1.} {\em If $n \ge 2$ and $({\mathcal I},|\ |_{\mathcal I})$ is a normed ideal so that ${\mathcal I} \ne {\mathcal C}_1$, then there are $n$-tuples of operators $\tau,\tau'$ so that $k_{\mathcal I}(\tau) = k_{\mathcal I}(\tau') = 0$ and $k_{\infty}^-(\tau \otimes \tau') \ne 0$.
}

\bigskip
\noindent
{\bf {\em Proof.}} It suffices to prove this for $n = 2$, since pairs of operators can be turned into $n$-tuples by repeating one of their components. Denoting by $\varphi(m)$ the number $|P|_{\mathcal I}$ where $P$ is an orthogonal projection of rank~$m$, the assumption ${\mathcal I} \ne {\mathcal C}_1$ is equivalent to $m^{-1}\varphi(m) \to 0$ as $m \to \infty$.

We construct two rooted trees with sets of vertices $X$ and $Y$ respectively, so that the order $w(v)$ of a vertex is either $2$ or $3$ and the order of the root is $2$ in both cases. We shall also denote by $d(v)$ the distance of a vertex $v$ to the root in both cases.

Let $(h_k)_{k \ge 1}$ be an increasing sequence of integers $\ge 1$ and let $S(p) = h_1 + \dots + h_p$.

The tree with vertices $X$ is such that $w(v) = 2$ if $S(2k-2) \le d(v) < S(2k-1)$, $k \ge 1$ and $w(v) = 3$ otherwise (except for the root). The tree with vertices $Y$ in turn is such that $w(u) = 2$ if $S(2k-1) \le d(u) < S(2k)$, $k \ge 1$ and $w(u) = 3$ otherwise (except for the root).

The canonical orthonormal bases of $\ell^2(X)$ and $\ell^2(Y)$ will be identified with $X$ and $Y$ respectively. On $\ell^2(X)$ we consider a pair of isometries $S_1,S_2$ so that $S_jv = v'$ where $d(v') = d(v) + 1$, $v \in X$ and $\{S_1v,S_2v\}$ is precisely the set of vertices adjacent to $v$ and which are one unit farther from the root than $v$. Note that $\{S_1v,S_2v\}$ has $2$ elements if $w(v) = 3$ and only $1$ element if $w(v) = 2$. We also consider a pair of isometries $T_1,T_2$ on $\ell^2(Y)$ satisfying the same kind of conditions.

The construction of the two trees is such that if $v \in X$, $u \in Y$ and $d(v) = d(u) \ge 1$ the $\{w(v),w(u)\} = \{2,3\}$. This implies that if $R_1 = S_1 \otimes T_1$, $R_2 = S_2 \otimes T_2$ and $\xi = v_0 \otimes u_0$, where $v_0,u_0$ are the roots, then 
\[
R_{\iota(1)}^{\alpha(1)} R_{\iota(2)}^{\alpha(2)} \dots R_{\iota(m)}^{\alpha(m)} \to R_{\iota(1)}^{\alpha(1)} R_{\iota(2)}^{\alpha(2)} \dots R_{\iota(m)}^{\alpha(m)} \xi
\]
where $m \ge 0$, $\iota(1) \ne \iota(2) \ne \dots \ne \iota(m)$, $\alpha(j) > 0$, $1 \le j \le m$ is an injection of the free monoid generated by $R_1,R_2$ into an orthonormal system of $\ell^2(X) \otimes \ell^2(Y)$. By (Prop.~3.4 in \cite{14}) this implies that $k_{\infty}^-(R_1,R_2) > 0$ and hence $k_{\infty}^-((S_1,T_1) \otimes (S_2,T_2)) > 0$.

On the other hand we will show that for a suitable choice of the sequence $(h_k)_{k \ge 1}$ also $k_{\mathcal I}(S_1,S_2) = k_{\mathcal I}(T_1,T_2) = 0$ will be satisfied.

Given $n \ge 1$ let $A(n)$ be the diagonal operator on $\ell^2(X)$ which is the multiplication operator by the function
\[
f(v) = \begin{cases}
1 &\mbox{if $d(v) \le S(2n-2)$} \\
(1 - (d(v) - S(2n-2))h_{2n-1}^{-1})_+ &\mbox{if $d(v) > S(2n-2)$.}
\end{cases}
\]
Remark that $f(v) = 0$ if $d(v) \ge S(2n-1)$. Thus $A(n) - S_j^*A(n)S_j$ is a diagonal operator with eigenvalues $0$ and $h_{2n-1}^{-1}$ and its rank is $\le 2^{S(2n-2)}h_{2n-1}$. It follows that
\[
|A(n) - S_j^*A(n)S_j|_{\mathcal I} \le h_{2n}^{-1}\varphi(2^{S(2n-2)}h_{2n-1}).
\]
Since $\varphi(m)m^{-1} \to 0$ as $m \to \infty$, when $h_1,\dots,h_{2n-2}$ are given we can find $h_{2n-1}$ large enough, so that
\[
h_{2n-1}^{-1}\varphi(2^{S(2n-2)}h_{2n-1}) \le n^{-1}.
\]
This insures that
\[
|A(n) - S_j^*A(n)S_j|_{\mathcal I} \to 0
\]
as $n \to \infty$ and since 
\[
|[A(n)S_j]|_{\mathcal I} = |S_j(A(n)-S_j^*A(n)S_j)|_{\mathcal I} \le |A(n)-S_j^*A(n)S_j|_{\mathcal I}
\]
and $A(n) \uparrow I$ as $n \to \infty$, we also get that
\[
k_{\mathcal I}(S_1,S_2) = 0.
\]

To satisfy $k_{\mathcal I}(T_1,T_2) = 0$ we consider $B(n)$ the diagonal operator on $\ell^2(Y)$ which is the multiplication operator by
\[
g(u) = \begin{cases}
1 &\mbox{if $d(u) \le S(2n-1)$} \\
1 - (d(u) - S(2n-1)h_{2n}^{-1})_+ &\mbox{if $d(u) > S(2n-1)$}
\end{cases}
\]
and remark that $g(u) = 0$ if $d(u) \ge S(2n)$.

Then $B(n) - T_j^*B(n)T_j$ is a diagonal operator with eigenvalues $0$ and $h_{2n}^{-1}$ and its rank is $\le 2^{S(2n-1)}h_{2n}$. We then choose $h_{2n}$ for given $h_1,\dots,h_{2n-1}$, so that
\[
h_{2n}^{-1}\varphi(2^{S(2n-1)}h_{2n}) \le n^{-1}.
\]
Like before for the $S_j$'s this will also insure that $k_{\mathcal I}(T_1,T_2) = 0$ is satisfied.\qed

\end{document}